\documentclass[a4paper,12pt]{article}

\usepackage{amsmath,amssymb}
\usepackage{amsfonts}

\normalfont

\def\qed{$\square$}
\newtheorem{theorem}{Theorem}[section]

\newtheorem{example}[theorem]{Example}
\newtheorem{lemma}[theorem]{Lemma}
\newtheorem{proposition}[theorem]{Proposition}

\newtheorem{remark}[theorem]{Remark}

\begin{document}

\title{Gaussian density estimates for  solutions to quasi-linear stochastic partial differential equations}
\author{
  \\
  {David Nualart}    \thanks{Supported by the NSF grant DMS-0604207.}         \\
  {\small \it Department of Mathematics }  \\[-0.15cm]
  {\small \it University of Kansas}          \\[-0.15cm]
  {\small \it Lawrence, Kansas, 66045, USA}  \\[-0.15cm]
  { \small {\tt nualart@math.ku.edu}}
  \\[-0.1cm]
\and
\\
  {Llu\'{\i}s Quer-Sardanyons \thanks{Supported by grants MEC-FEDER Ref. MTM20006-06427 from the
  Direcci\'on General de Investigaci\'on, Ministerio de Educaci\'on y Ciencia, Spain, and BE 2007 from the Ag\`encia
  de Gesti\'o
  d'Ajuts Universitaris i de Recerca, Generalitat de Catalunya.}}              \\
  {\small\it Departament de Matem\`atiques} \\[-0.15cm]
  {\small\it  Universitat Aut\`onoma de Barcelona}          \\[-0.15cm]
  {\small\it 08193 Bellaterra (Barcelona), Spain}  \\[-0.15cm]
  {\small  {\tt quer@mat.uab.cat }} \\[-0.1cm]
  \\}

\maketitle

\begin{abstract}
In this paper we establish lower and upper Gaussian bounds for the
solutions to the heat and wave equations driven by an additive
Gaussian noise, using the techniques of Malliavin calculus and
recent density estimates obtained by Nourdin and Viens in \cite{nv}.
In particular, we deal with the one-dimensional stochastic heat
equation in $[0,1]$ driven by the space-time white noise, and the
stochastic heat and wave equations in $\mathbb{R}^d$ ($d\geq 1$ and
$d\leq 3$, respectively) driven by a Gaussian noise which is white
in time and has a general spatially homogeneous correlation.
\end{abstract}

\section{Introduction}
The Malliavin calculus, or stochastic calculus of variations, is a
suitable technique for proving that a given random vector $F=(F_1,
\dots, F_m)$ on the Wiener space possesses a smooth  probability
density. There has been a current interest in the applications of
the stochastic  calculus of variations to obtain lower and upper
bounds of Gaussian type for the density of a given Wiener
functional.  The starting point of this research is the work by
Kusuoka and Stroock \cite{ks}, where they proved that the density of
a uniformly hypoelliptic diffusion whose drift is a smooth
combination of its diffusion coefficient has a lower bound of
Gaussian type. Recently, three different approaches have been
developed to derive  Gaussian-type bounds for  densities of general
Wiener functionals using the techniques of Malliavin calculus:
\begin{itemize}
\item[(i)] In the paper \cite{k}, Kohatsu-Higa proposes  a general
methodology for computing lower bounds for a multidimensional
functional $F=(F_1, \dots, F_m)$ of   a  Wiener sheet in
$[0,T]\times \mathbb{R}^d$.  The method uses an approximation of $F$
by means of a sequence of conditionally non-degenerate random
variables adapted to the filtration generated by the white noise.
This paper was inspired by the work by Kusuoka and Stroock
\cite{ks}, and it can be used  in  a non Markovian  framework. As an
application, the author obtains lower bounds for the solution to the
one-dimensional heat equation driven by a space-time white noise,
assuming that the diffusion coefficient $\sigma$ is bounded away
from zero. The ideas introduced in this paper were later developed
in the work by Bally \cite{bally} to get Gaussian  lower bounds for
locally elliptic It\^o processes.
\item[(ii)]
 Nourdin and Viens in \cite{nv} have proved a new formula for the density of a one-dimensional
 Wiener functional in terms of the Malliavin calculus (see (\ref{22}) below). As an application, they obtain upper and lower Gaussian bounds for the density of the maximum of a general Gaussian process.
 \item[(iii)]
  In  \cite{mn} Malliavin and Nualart      derived Gaussian lower bounds for multidimensional Wiener functionals under an exponential moment condition on the divergence of a covering vector field. A one-dimensional version of this result was obtained by Nualart in \cite{enu}.
\end{itemize}

The purpose of this paper is to apply the results obtained by
Nourdin and Viens in \cite{nv} to the solutions of different classes
of stochastic partial differential equations driven by an additive
Gaussian noise.  Upper and lower Gaussian estimates for the density
of solutions to stochastic partial differential equations are
essential tools  in the potential analysis for these type of
equations (see the recent works \cite{dkn1, dkn2, dn}).

The paper is organized as follows.
 In Section 2 we recall briefly the results of \cite{nv}.
Section 3 deals with the one-dimensional heat equation on the
interval $[0,1]$ with Dirichlet boundary conditions and driven by a
space-time white noise. Let us denote by $u(t,x)$ its solution
evaluated at some $(t,x)\in \mathbb{R}_+\times (0,1)$ and set
$m=E|u(t,x)|$. Then, we derive a gaussian lower and upper bound for
the density  of $u(t,x)$ of the form
\begin{equation}  \label{eq1}
 Ct^{-\frac  12} \exp \left( - \frac {(z- m)^2}{C' t^{\frac 12}}\right).
\end{equation}
The lower bound differs from that obtained by Kohatsu-Hida in the
paper \cite{k},  where the  first factor was $Ct^{\frac 14}$.

Section 3 is devoted to the stochastic heat equation in
$\mathbb{R}^d$ with an additive Gaussian  noise which is white in
time and it has a spatially homogeneous correlation in space. We
assume that the spectral measure of the noise  integrates $(1+
|\xi|^2)^\eta$ for some $\eta\in (0,1)$, and we obtain  lower and
upper Gaussian bounds for the density of the solution of the form
(\ref{eq1}), that involve the powers $t $ and $t^{1-\eta}$ (see
Theorem \ref{bounds} for the precise statement). Let us mention at
this point that it has been of much importance the fact the
Malliavin derivative of the solution defines a non-negative process,
since it solves a linear parabolic equation.

Finally, in Section 4 we consider a stochastic wave equation in
dimension $d=1,2,3$, again driven by an additive  Gaussian noise
which is white in time and it has a spatially homogeneous
correlation in space. In this case we obtain  lower and upper bounds
only for $t$ small enough, and they involve the powers $t^3 $ and
$t^{3-2\eta}$, where again $\eta$ is given by the integrability of
the spectral measure. We have not been able to overcome this
technical restriction on the time parameter because, in comparison
with the heat equation, the Malliavin derivative of the solution in
this case does not need to be a non-negative function neither a
function itself.

\section{Gaussian density estimates}
\label{gaussian}

In this section, we   recall a general method set up in \cite{nv}
in order to show that a smooth   random
variable in the sense of Malliavin calculus has a probability density admitting Gaussian lower and upper bounds.

First of all, let us briefly describe the Gaussian context in which
we will be working on and the main elements of the Malliavin
calculus associated to it. Namely, suppose that in a complete
probability space $(\Omega,\mathcal{F},P)$ we are given a centered
Gaussian family $W=\{W(h),\; h\in H\}$ of random variables  indexed by a separable
Hilbert space $H$ with covariance
$$E(W(h)W(g))=\langle h,g\rangle_H,\; h,g\in H.$$
The family $W$ is usually called an isonormal Gaussian process on
$H$. Assume that $\mathcal{F}$ is the $\sigma$-field generated by
$W$.

Let us use standard notation for the main operators of the Malliavin
calculus determined by the family $W$ (see, for instance,
\cite{nualart}). More precisely, we denote by $D$ the Malliavin
derivative, defined as a closed and unbounded operator from
$L^2(\Omega)$ into $L^2(\Omega; H)$, whose domain is denoted by
$\mathbb{D}^{1,2}$. The adjoint of the operator $D$ is denoted by
$\delta$, usually called the divergence operator. A random element
$u\in L^2(\Omega;H)$ belongs to the domain of $\delta$ if and only
if it satisfies
$$|E(\langle DF,\delta(u)\rangle_H)|\leq C \|F\|_{L^2(\Omega)},$$
for any $F\in \mathbb{D}^{1,2}$, where the constant $C$ only depends
on $u$. For any element $u$ in the domain of $\delta$, the random
variable $\delta(u)$ can be characterized by the duality
relationship
$$E(F \delta(u))=E(\langle DF,u\rangle_H),$$
for every $F\in \mathbb{D}^{1,2}$.

Any random variable $F$ in $L^2(\Omega,\mathcal{F},P)$ can be
decomposed by means of its Wiener chaos expansion (see \cite[Section
1.1.1]{nualart}), which is usually written as
$$F=\sum_{n=0}^\infty J_n F,$$
where $J_n$ denotes the projection onto the $n$th Wiener chaos.

Using the chaos expansion one may define the operator $L$ by the
formula $L=\sum_{n=0}^\infty -n J_n$, which is called the generator
of the Orstein-Uhlenbeck semigroup (see \cite[Section
1.4]{nualart}). It is related to the Malliavin derivative $D$ and
its adjoint $\delta$ through the formula
$$\delta DF=-LF,$$
in the sense that $F$ belongs to the domain of $L$ if and only if it
belongs to the domain of $\delta D$, and in this case the above
equality holds.

One can also define the inverse of $L$, denoted by $L^{-1}$, as
follows: for any $F\in L^2(\Omega,\mathcal{F},P)$, set $L^{-1}
F:=\sum_{n=1}^\infty -\frac{1}{n} J_n F$. Then it holds that $L
L^{-1}F=F-E(F)$, for any $F\in L^2(\Omega,\mathcal{F},P)$, so that
$L^{-1}$ acts as the inverse of $L$ for centered random variables.\\

Let us consider $F\in \mathbb{D}^{1,2}$ with mean zero and we define
the following function in $\mathbb{R}$:
$$g(z):=E [ \langle DF,-DL^{-1} F\rangle_H |F=z].$$
By \cite[Proposition 3.9]{np}, it holds that $g(z)\geq 0$ on the
support of $F$. Then, \cite[Theorem 3.1]{nv} states that if the
random variable $g(F)$ is bounded away from zero almost surely, that
is
$$g(F)\geq c_1>0,\; \text{a.s.},$$
for some constant $c_1$,
then $F$ has a density $\rho$ whose support is $\mathbb{R}$
satisfying, almost everywhere: \begin{equation}
\rho(z)=\frac{E|F|}{2g(z)}\exp\left(-\int_0^z
\frac{y}{g(y)}dy\right). \label{22} \end{equation} As a consequence,
see \cite[Corollary 3.4]{nv}, if one also has that $g(F)\leq c_2$,
a.s., then the density $\rho$ satisfies, for almost all $z\in
\mathbb{R}$:
\begin{equation}
\label{eq2}
\frac{E|F|}{2c_1}\exp\left(-\frac{z^2}{2c_2}\right)\leq
\rho(z)\leq \frac{E|F|}{2c_2}\exp\left(-\frac{z^2}{2c_1}\right).
\end{equation}

Let us also mention that \cite[Proposition 3.5]{nv}
provides a more suitable formula for $g(F)$ for computational
purposes. Indeed, given a random variable $F\in \mathbb{D}^{1,2}$,
one can write $DF=\Phi_F(W)$, where $\Phi_F$ is a measurable mapping
from $\mathbb{R}^H$ to $H$, determined $P\circ W^{-1}$-almost surely
(see \cite{nualart}, p. 54-55). Then \begin{equation}
g(F)=\int_0^\infty e^{-\theta} \textbf{E}\left[ \langle
\Phi_F(W),\Phi_F(e^{-\theta}W+\sqrt{1-e^{-2\theta}}W')\rangle_H
\big| F\right] d\theta, \label{23} \end{equation} where $W'$ stands
for an independent copy of $W$ such that $W$ and $W'$ are defined on
the product probability space $(\Omega\times
\Omega',\mathcal{F}\otimes \mathcal{F}',P\times P')$. Eventually,
$\textbf{E}$ denotes the mathematical expectation with respect to
$P\times P'$.

Formula (\ref{23}) can be still rewritten in the following form:
\begin{equation}
g(F)= \int_0^\infty  e^{-\theta} E\left[ E'\left( \left\langle
DF,\widetilde{DF}\right\rangle_H\right) \big| F \right] d\theta,
\label{8}
\end{equation}
where, for any random variable $X$ defined in
$(\Omega,\mathcal{F},P)$, $\widetilde{X}$ denotes the shifted random
variable in $\Omega\times \Omega'$
$$\widetilde{X}(\omega,\omega')=X(e^{-\theta}\omega+
\sqrt{1-e^{-2\theta}}\omega'),\; \omega\in \Omega,\; \omega'\in
\Omega.$$ Notice that, indeed, $\widetilde{X}$ depends on the
parameter $\theta$, but we have decided to drop its explicit
dependence for the sake of simplicity.

Along the paper we denote by $C$ a generic constant which may  vary from line to line.

\section{Density estimates for the stochastic heat equation in $[0,1]$}
\label{boundary}

In this section, we will consider a stochastic heat equation in
$[0,1]$, with Dirichlet boundary conditions, some non-linear drift
$b$ and with an additive space-time white noise perturbation. We aim
to give sufficient conditions on the coefficient $b$ ensuring
Gaussian lower and upper bounds for the probability density of the
solution at any point.

\subsection{The stochastic heat equation in $[0,1]$}
\label{setting-boundary}

We are concerned with the following one dimensional heat equation
driven by a space-time white noise:
\begin{equation}
\frac{\partial u}{\partial t}(t,x)-\frac{\partial ^2 u}{\partial
x^2} (t,x)=b(u(t,x)) + \sigma \dot W(t,x),\quad (t,x)\in [0,T]\times
[0,1], \label{1}
\end{equation}
where $T>0$, the initial condition is given by a continuous function
$u_0:[0,1]\rightarrow \mathbb{R}$ and we consider Dirichlet boundary
conditions. That is,
\begin{equation}
u(0,x)=u_0(x),\; x\in [0,1], \label{2}
\end{equation}
$$u(t,0)=u(t,1)=0,\; t\in [0,T].$$
The real-valued random field solution to Equation (\ref{1}) will be
denoted by $\{u(t,x),\;  (t,x)\in [0,T]\times[0,1]\}$. The function
$b:\mathbb{R}\rightarrow \mathbb{R}$ is of class $\mathcal{C}^1$
having a bounded derivative and $\sigma> 0$ is a constant. We assume that
$\{W(t,x), \; (t,x)\in [0,T]\times [0,1]\}$ is a Brownian sheet
on $[0,T]\times [0,1]$, defined in a complete probability space $(\Omega, \mathcal{F}, P)$.  That is, $\{W(t,x)\}$ is a centered Gaussian family with the covariance function
\[
E(W(t,x)W(s,y))= (t\wedge s)(x\wedge y).
\]
For $0\leq t\leq T$, let $\mathcal{F}_t$ be
the $\sigma$-field generated by the random variables $\{W(s,x), \;
(s,x)\in [0,t]\times [0,1]\}$  and  the $P$-null sets.

The solution to the formal Equation (\ref{1}) is understood in the
mild sense: a $\{\mathcal{F}_t\}$-adapted stochastic process
$\{u(t,x),\; (t,x)\in [0,T]\times [0,1]\}$ solves (\ref{1}) with
initial and boundary conditions (\ref{2}) if, for any $(t,x)\in
(0,T]\times (0,1)$,
\begin{align}
u(t,x) = & \int_0^1 G_t(x,y) u_0(y) dy+ \int_0^t \int_0^1
G_{t-s}(x,y) b(u(s,y)) dy
ds\nonumber \\
& \quad + \sigma \int_0^t\int_0^1 G_{t-s}(x,y) W(ds,dy), \label{3}
\end{align}
where $G_t(x,y)$, $(t,x,y)\in \mathbb{R}_+\times (0,1)^2$, denotes
the Green function associated to the heat equation on $[0,1]$ with Dirichlet
boundary conditions.  We will use the following facts,   for any $0<t'<t $ and $x\in (0,1)$,
\begin{equation}
0\leq G_t(x,y)\leq  \frac 1{ \sqrt{4\pi t} }   e^{-\frac{(x-y)^2}{4t}},  \label{4}
\end{equation}
\begin{equation}
c_x \sqrt{t-t'} \le \int_{t'}^t \int_0^1 |G_{t-s}(x,y)|^2 dy ds \leq
\frac{1}{\sqrt{2\pi}} \sqrt{t-t'},  \label{5}
\end{equation}
where $c_x$ is a positive constant.  We also have $\inf_{\alpha \le
x \le 1-\alpha}c_x>0$ for any $\alpha \in (0,1)$. In the case of
Neumann boundary conditions,  $c_x$ does not depend on $x$ (see, for
instance,
  \cite[Lemma A1.2]{morien}).

Let us also mention that the stochastic integral in (\ref{3}) is
understood as an integral with respect to the Brownian sheet in the
sense of Walsh \cite{walsh}.

Existence and uniqueness of mild solution for Equation (\ref{3}) can
be deduced from the results in \cite{walsh}; for an even more
general setting, see also \cite{bgp}. Malliavin calculus applied to
Equation (\ref{3}) have been dealt with in \cite{bp}. In this case,
we consider the  Gaussian context associated to the space-time white
noise. That is, we have $H=L^2(\mathbb{R}_+\times [0,1])$ and the
Gaussian family $\{W(h),\; h\in L^2(\mathbb{R}_+\times [0,1])\}$ is
given by the Wiener integral
$$W(h)=\int_{\mathbb{R}_+} \int_0^1 h(s,y) W(ds,dy).$$

A consequence of Proposition 4.3 and Theorem 2.2 in \cite{bp} is
that, for all $(t,x)\in (0,T]\times [0,1]$, the random variable
$u(t,x)$ belongs to $\mathbb{D}^{1,2}$, the Malliavin derivative
satisfies the linear parabolic equation
\begin{equation}
D_{r,z}u(t,x)= \sigma G_{t-r}(x,z) + \int_r^t \int_0^1 G_{t-s}(x,y)
b'(u(s,y)) D_{r,z}u(s,y)\; dy ds, \label{7}
\end{equation}
for any $(r,z)\in [0,t]\times [0,1]$ and $(t,x)\in [0,T]\times
[0,1]$, and the probability law of $u(t,x)$ has a density. The
positivity of $\sigma$ and $G$ guarantees that the solution of
Equation (\ref{7}) remains non-negative, that is $D_{r,z}u(t,x)\geq
0$, a.s. This will be a key point in proving that the
density of $u(t,x)$ admits Gaussian upper and lower estimates.


\subsection{Gaussian bounds for the density} \label{main}

We fix $T>0$ and consider $\{u(t,x),\; (t,x)\in [0,T]\times [0,1]\}$
the unique mild solution of Equation (\ref{3}). This section will be
devoted to proof the following result.

\begin{theorem}  \label{teo1}
Assume that the drift coefficient $b$ is a $\mathcal{C}^1$ function
with bounded derivative. Then, for all $t\in (0,T]$ and $x\in
(0,1)$, the random variable $u(t,x)$ possesses a density $p$
satisfying the following statement: for almost every $z\in
\mathbb{R}$,
\begin{equation} \frac{E\left| u(t,x)-m \right|}{ C_2 t^{\frac{1}{2}}}
\exp\left\{ -\frac{(z-m)^2}{C_1 t^{\frac{1}{2}} }\right\}\leq p(z)
\leq \frac{E\left| u(t,x)-m\right|}{ C_1 t^{\frac{1}{2}}}
\exp\left\{ -\frac{(z-m)^2}{ C_2 t^{\frac{1}{2}} }\right\},
\label{9}
\end{equation} where $m:=E(u(t,x))$ and $C_1$, $C_2$
are positive quantities  depending on $\sigma$, $\|b'\|_\infty$, $T$ and $x$.
\end{theorem}

The statement of the above Theorem \ref{main} will be a consequence
of formula (\ref{22}) (see \cite[Theorem 3.1]{nv}) and the following
proposition. In order not to overload notations, set
$F:=u(t,x)-E(u(t,x))$, and recall that, as it has been specified in
Section \ref{gaussian} (see (\ref{8}) therein),
\begin{align}
g(F) &= \int_0^\infty  \; e^{-\theta} E\left[ E'\left( \int_0^t
\int_0^1 D_{r,z}F \left(\widetilde{D_{r,z}F}\right) dz dr
\right)\Big| F\right] d\theta \nonumber \\
& =\int_0^\infty \; e^{-\theta} E\left[ E'\left( \int_0^t \int_0^1
D_{r,z}u(t,x) \left(\widetilde{D_{r,z}u(t,x)}\right) dz dr
\right)\Big| F\right] d\theta, \label{13}
\end{align} where $\widetilde{DF}=(DF)(e^{-\theta}\omega +
\sqrt{1-e^{-2\theta}}\omega')$.

\begin{proposition}\label{prop}
Fix $T>0$ and assume that the drift coefficient $b$ is of class
$\mathcal{C}^1$ and has a bounded derivative. There exist positive
constants $C_1, C_2$ such that
\begin{equation}
C_1 t^{\frac{1}{2}} \leq g(F) \leq C_2 t^{\frac{1}{2}}, \label{14}
\end{equation}
for all $t\in (0,T]$.
\end{proposition}

In order to prove Proposition \ref{prop}, we will need the following
technical lemma:

\begin{lemma}
Let $t>0$. Assume that $b\in \mathcal{C}^1$ and it has a bounded
derivative. There exists a positive constant $K$  depending on $\sigma$ and $\|b'\|_\infty$, such that, for any
$\delta\in (0,1]$: \begin{equation} \sup_{\stackrel{(1-\delta)t\leq
\nu\leq t}{0\leq y \leq 1}} \int_{(1-\delta)t}^t \int_0^1 E\left[
\left| D_{r,z}u(\nu,y)\right|^2 \big| F\right] \; dz dr\leq K
(\delta t)^{\frac 1 2}\label{10}\end{equation} and
\begin{equation}\sup_{\theta\in \mathbb{R}} \sup_{\stackrel{(1-\delta)t\leq \nu\leq
t}{0\leq y \leq 1}} \int_{(1-\delta)t}^t \int_0^1 E\left[ E'\left(
\left| \widetilde{D_{r,z}u(\nu ,y)}\right|^2 \right) \Big| F\right]
\; dz dr\leq K (\delta t)^{\frac 1 2},\label{11}\end{equation}
$P$-almost surely. \label{lemma1}
\end{lemma}
\emph{Proof:} We will only deal with the proof of (\ref{10}), since
(\ref{11}) may be checked using exactly the same arguments.

Let us first invoke the linear equation (\ref{7}) satisfied by the
Malliavin derivative $Du(\nu,v)$, for $(\nu,v)\in
[(1-\delta)t,t]\times (0,1)$, and then take the square $L^2$-norm on
$[(1-\delta)t,t]\times [0,1]$:
\begin{align}
& \int_{(1-\delta)t}^t \int_0^1 |D_{r,z}u(\nu,v)|^2\; dz dr  \leq 2\sigma^2
\int_{(1-\delta)t}^t \int_0^1  |G_{\nu-r}(v,z)|^2 \; dz dr \nonumber
\\
& \quad \quad \quad + 2 \int_{(1-\delta)t}^t \int_0^1 \left(
\int_r^\nu \int_0^1 G_{\nu-s}(v,y) b'(u(s,y)) D_{r,z}u(s,y) \; dy
ds\right)^2\; dz dr. \label{12}
\end{align}
By (\ref{5}), the first term in the right-hand side of (\ref{12})
can be bounded by $\sqrt{ \frac 2\pi} \, \sigma^2 (\delta t)^{\frac
1 2}$. For the second one, we apply H\"older's inequality, the fact
that $b'$ is bounded and Fubini's theorem, so that we end up with
\begin{align*}
&  \int_{(1-\delta)t}^t \int_0^1 |D_{r,z}u(\nu,v)|^2\; dz dr \leq  \sqrt{ \frac 2\pi} \, \sigma^2 (\delta t)^{\frac 1 2}   \\
&  \qquad  +  2 \|b'\|_\infty^2 \delta t \int_{(1-\delta)t}^\nu
\int_0^1 |G_{\nu-s}(v,y)|^2  \left( \int_{(1-\delta)t}^t \int_0^1
|D_{r,z}u(s,y)|^2\; dz dr\right) \; dy ds.
\end{align*}
Taking the conditional expectation $E[\;\cdot\;|F]$ and using again
(\ref{5}) we obtain
\begin{align*}
&  \sup_{\stackrel{(1-\delta)t\leq \rho \leq \nu}{0\leq v \leq 1}}
\int_{(1-\delta)t}^t \int_0^1 E\left[ |D_{r,z}u(\rho,v)|^2 \big|
F\right] \; dz dr \leq \sqrt{ \frac 2\pi} \, \sigma^2  (\delta
t)^{\frac 1 2}
\\
& \quad \quad  +  2 \|b'\|_\infty^2  \delta t \int_{(1-\delta)t}^\nu  \left(
\sup_{\stackrel{(1-\delta)t\leq \tau\leq s}{0\leq y \leq 1}}
\int_{(1-\delta)t}^t \int_0^1 E\left[ |D_{r,z}u(\tau,y)|^2 \big|
F\right] \; dz dr\right) \frac{1}{\sqrt{\nu-s}}\; ds.
\end{align*}
If we set
$$\Psi_{\delta,t}(\nu):=\sup_{\stackrel{(1-\delta)t\leq \rho \leq \nu}{0\leq v \leq 1}}
\int_{(1-\delta)t}^t \int_0^1 E\left[ |D_{r,z}u(\rho,v)|^2 \big|
F\right] \; dz dr ,\; \nu\in [(1-\delta)t,t],$$ then we have seen
that
$$\Psi_{\delta,t}(\nu)\leq K\, (\delta t)^{\frac 1 2} + K\, \delta t
\int_{(1-\delta)t}^\nu \Psi_{\delta,t}(s) \frac{1}{\sqrt{\nu-s}}\;
ds, \; \nu\in [(1-\delta)t,t],\, a.s.$$ Now we can conclude by
applying Gronwall's lemma \cite[Lemma 15]{Da}. \qed

\vspace{0.4cm}

\emph{Proof of Proposition \ref{prop}:} We first recall that the
Malliavin derivative of $u(\nu,v)$, $(\nu,v)\in [0,T]\times [0,1]$,
satisfies that $D_{r,z}u(\nu,v)\geq 0$, for all $(r,z)\in
[0,T]\times [0,1]$, a.s. This is because the Malliavin derivative
solves the linear parabolic equation (\ref{7}). Let us deal with the
proof of
(\ref{14}) in two steps:\\

\emph{Step 1: The lower bound.} Fix $\delta \in (0,1]$ and let us
first derive the lower bound in (\ref{14}). Since the Malliavin
derivative of $u(t,x)$ is non-negative, formula (\ref{13}) yields
$$g(F)\geq \int_0^\infty  \; e^{-\theta} E\left[
E'\left( \int_{(1-\delta)t}^t \int_0^1 D_{r,z} u(t,x)
\left(\widetilde{D_{r,z} u(t,x)}\right)\; dz dr \right)\Big|
F\right]d\theta.$$ By Equation (\ref{7}), we can decompose the
right-hand side of the above inequality in a sum of four terms:
$$A_0(t,x;\delta)=\sigma^2 \int_{(1-\delta)t}^t \int_0^1
|G_{t-r}(x,z)|^2 \; dz dr,$$
\begin{align*}
A_1(t,x;\delta)= & \sigma \int_{(1-\delta)t}^t \int_0^1
G_{t-r}(x,z)\\ & \quad \quad \times E\left[ \int_r^t \int_0^1
G_{t-s}(x,y) b'(u(s,y)) D_{r,z}u(s,y) \; dy ds \; \Big| F\right]\;
dz dr,
\end{align*}
\begin{align*}
A_2(t,x;\delta)= & \sigma \int_0^\infty e^{-\theta}
\int_{(1-\delta)t}^t \int_0^1 G_{t-r}(x,z)\\ & \; \times E\left[
E'\left( \int_r^t \int_0^1 G_{t-s}(x,y) b'(\widetilde{u(s,y)})
\left(\widetilde{D_{r,z}u(s,y)}\right) \; dy ds \right)\; \Big|
F\right]\; dz dr d\theta,
\end{align*}
\begin{align*}
A_3(t,x;\delta)= & \int_0^\infty e^{-\theta} \int_{(1-\delta)t}^t
\int_0^1 E\left[ \left(\int_r^t \int_0^1 G_{t-s}(x,y) b'(u(s,y))
D_{r,z}u(s,y) \; dy ds \right) \right.\\
& \; \times \left. E'\left( \int_r^t \int_0^1 G_{t-s}(x,y)
b'(\widetilde{u(s,y)}) \left(\widetilde{D_{r,z}u(s,y)}\right) \; dy
ds \right)\; \Big| F\right]\; dz dr d\theta.
\end{align*}
First we notice that, by (\ref{5}): \begin{equation}
A_0(t,x;\delta)\geq \sigma^2 c_x (\delta t)^{\frac{1}{2}},
\label{15}
\end{equation}
 Thus we can write \begin{equation}
g(F)\geq  \sigma^2 c_x (\delta t)^{\frac{1}{2}}
-|A_1(t,x;\delta)+A_2(t,x;\delta)+A_3(t,x;\delta)|, \label{16}
\end{equation} so that we will need to obtain upper bounds for the
terms
$|A_i(t,x;\delta)|$, $i=1,2,3$.\\

We apply Fubini's theorem for the conditional expectation, the
boundedness of $b'$, Cauchy-Schwarz's inequality and the bound
(\ref{5}), so we have the following estimates:
\begin{align*}
& |A_1(t,x;\delta)|\nonumber \\
& \; \leq  C \left( \int_{(1-\delta)t}^t \int_0^1
|G_{t-r}(x,z)|^2\; dz dr\right)^{\frac{1}{2}}\nonumber  \\
&\quad  \quad \quad   \times \left( \int_{(1-\delta)t}^t \int_0^1
\left| \int_r^t \int_0^1 G_{t-s}(x,y) E\left[ |D_{r,z}u(s,u)|\;\big|
F\right]
dyds\right|^2 dz dr\right)^{\frac{1}{2}}\nonumber \\
& \; \leq C (\delta t)^{\frac{3}{4}} \left( \int_{(1-\delta)t}^t
\int_0^1 \left( \int_r^t \int_0^1 |G_{t-s}(x,y)|^2 E\left[
|D_{r,z}u(s,y)|^2 \; \big| F\right]  dy ds
\right) dz dr \right)^{\frac{1}{2}}\nonumber \\
& \; \leq C (\delta t)^{\frac{3}{4}}  \left( \int_{(1-\delta)t}^t
\int_0^1 |G_{t-s}(x,y)|^2 \left( \int_{(1-\delta)t}^t \int_0^1
E\left[ |D_{r,z}u(s,y)|^2 \; \big| F\right]  dz dr \right) dy ds
\right)^{\frac{1}{2}} \nonumber \\
& \; \leq C (\delta t) \left( \sup_{\stackrel{(1-\delta)t\leq s\leq
t}{0\leq y \leq 1}} \int_{(1-\delta)t}^t \int_0^1 E\left[
|D_{r,z}u(s,y)|^2 \; \big| F\right]  dz dr \right)^{\frac{1}{2}}.
\end{align*}
At this point we are in position to apply (\ref{10}) in Lemma
\ref{lemma1}. Therefore \begin{equation} |A_1(t,x;\delta)|\leq C\;
(\delta t)^{\frac{5}{4}},\; a.s., \label{17} \end{equation} where
$C=\frac{1}{\sqrt{\pi} (2\pi)^{\frac 14}} \, \sigma^2 \|b'\|_\infty $.\\

In order to get a bound for $|A_2(t,x;\delta)|$, one can use
analogous arguments as for $|A_1(t,x;\delta)|$ but applying
(\ref{11}) instead of (\ref{10}) in Lemma (\ref{lemma1}). Hence, one
obtains \begin{equation} |A_2(t,x;\delta)|\leq C\; (\delta
t)^{\frac{5}{4}},\; a.s. \label{18} \end{equation}

\vspace{0.3cm}

Let us finally estimate $|A_3(t,x;\delta)|$. For this, we apply
Fubini's theorem, the fact that $b'$ is bounded, Cauchy-Schwarz
inequality with respect to $dz dr dP_{|F} dP'$ and we finally invoke
Lemma \ref{lemma1}: \begin{align}
& |A_3(t,x;\delta)|\nonumber  \\
& \; \leq C\int_0^\infty e^{-\theta} \int_{(1-\delta)t}^t \int_0^1
\int_{(1-\delta)t}^t \int_0^1
G_{t-s}(x,y) G_{t-\bar s}(x,\bar y)\nonumber \\
& \quad \times  \left( \int_{(1-\delta)t}^t\int_0^1 E\left[ E'\left(
|D_{r,z}u(s,y) \widetilde{D_{r,z}u(\bar s,\bar y)}| \right) \Big|
F\right] dz
dr\right) dy ds d\bar y d\bar s d\theta \nonumber \\
& \; \leq C\int_0^\infty e^{-\theta} \int_{(1-\delta)t}^t \int_0^1
\int_{(1-\delta)t}^t \int_0^1
G_{t-s}(x,y) G_{t-\bar s}(x,\bar y)\nonumber \\
& \quad \times  \left( \int_{(1-\delta)t}^t\int_0^1 E\left[
|D_{r,z}u(s,y)|^2  \Big| F\right] dz
dr\right)^{\frac{1}{2}}\nonumber \\
& \quad \times \left( \int_{(1-\delta)t}^t\int_0^1 E\left[ E'\left(
|\widetilde{D_{r,z}u(\bar s,\bar y)}|^2 \right) \Big| F\right]
\right)^{\frac{1}{2}} dy ds d\bar y d\bar s d\theta \nonumber \\
& \; \leq  C (\delta t)^{\frac{1}{2}} \left( \int_{(1-\delta)t}^t
\int_0^1
G_{t-s}(x,y) dy ds\right)^2 \nonumber \\
& \; \leq C (\delta t)^2, \label{19}
\end{align}
where $C=\frac{1}{\sqrt{2} \pi^{\frac 32}} \, \sigma^2
\|b'\|_\infty^2$. The very last estimate in (\ref{19}) has been
obtained after applying Cauchy-Schwarz's inequality and the bound
(\ref{5}). Eventually, plugging the bounds (\ref{17})-(\ref{19}) in
(\ref{16}) we have
$$
g(F)\geq  \sigma^2 c_x(\delta t)^{\frac{1}{2}}- c_1 \left((\delta
t)^{\frac{5}{4}} + (\delta t)^2\right),
$$
where $c_1$ is a positive constant depending on $\sigma$ and $\|b'\|_\infty$.
Hence, if
we assume that $\delta< 1\wedge \frac{1}{T}$, then we can write
$$
g(F)\geq \sigma^2 c_x(\delta t)^{\frac{1}{2}}- 2 c_1 (\delta
t)^{\frac{5}{4}} \geq \left( \sigma^2 c_x \delta^{\frac{1}{2}} - 2 c_1
\delta^{\frac{5}{4}} T^{\frac{3}{4}} \right) t^{\frac{1}{2}}.
$$ It
only remains to observe that the quantity $\sigma^2 c_x \delta^{\frac{1}{2}}
- 2 c_1 \delta^{\frac{5}{4}} T^{\frac{3}{4}}$ is strictly positive
whenever $\delta $ is sufficiently small, namely $\delta \in
(0,\delta_0)$, with
$$\delta_0=1\wedge \frac{1}{T}\wedge \frac{1}{T}
\left( \frac{\sigma^2 c_x}{2c_{1}}\right)^{\frac{4}{3}}.$$ Thus,
the lower bound in (\ref{14}) has been then proved.\\

\emph{Step 2: The upper bound.} The upper estimation in (\ref{14})
is almost an immediate consequence of the computations which we have
just performed for the lower bound. More precisely, according to
(\ref{13}) and the considerations in the first part of the proof, we
have the following: \begin{equation} g(F)\leq \sum_{i=0}^3
|A_i(t,x;1)|,\label{20}\end{equation} where we notice that we have
substituted $\delta$ by $1$ in $A_i(t,x;\delta)$, $i=0,1,2,3$. We
have already seen that $|A_i(t,x;1)|\leq C t^{\frac{5}{4}}$, for
$i=1,2$, and $|A_3(t,x;\delta)|\leq C t^2$, so we just need to bound
$|A_0(t,x;1)|$, which follows directly from (\ref{5}). Thus
$$g(F)\leq C\left( t^{\frac{1}{2}}+t^{\frac{5}{4}} + t^2\right),\; a.s.,$$
for a constant $C>0$. Therefore $g(F) \leq C_2 t^{\frac{1}{2}}$,
where the constant $C_2>0$   depends on $\sigma$, $\|b'\|_\infty$
and $T$.
 \qed

\vspace{0.3cm}

We are now in position to prove the main result of this section:

\vspace{0,3cm}

\emph{Proof of Theorem \ref{main}:} The random variable
$F=u(t,x)-E(u(t,x))$ is centered, belongs to $\mathbb{D}^{1,2}$ and,
by Proposition \ref{prop}, it holds that $0<C_1 t^{\frac 1 2} \leq
g(F)$, for all $t\in (0,T]$. We apply then \cite[Theorem 3.1]{nv}
and we obtain that the probability density
$\rho:\mathbb{R}\rightarrow \mathbb{R}$ of $F$ is given by
$$ \rho(z)=\frac{E|u(t,x)-E(u(t,x))|}{2g(z)}\exp\left(-\int_0^z
\frac{y}{g(y)}dy\right),$$ for almost every $z\in \mathbb{R}$. Thus,
the density $p$ of the random variable $u(t,x)$ satisfies
\begin{equation} p(z)=
\frac{E|u(t,x)-E(u(t,x))|}{2g(z-E(u(t,x)))}\exp\left(-\int_0^{z-E(u(t,x))}
\frac{y}{g(y)}dy\right). \label{21} \end{equation} In order to
conclude the proof, we only need to use the bounds (\ref{14}) in the
above expression (\ref{21}). \qed


\section{The stochastic heat equation in $\mathbb{R}^d$}
\label{initial}

In this section  we are interested in the stochastic heat equation
in $\mathbb{R}^d$ with an additive Gaussian noise which is white in
time and it has a spatially homogeneous correlation in space. We aim
to find out sufficient conditions on the drift term and the noise's
spatial correlation ensuring that the density of the solution admits
Gaussian-type estimations. The fact that we deal with a SPDE in
$\mathbb{R}^d$ with a non-trivial spatial correlation makes the
analysis in this case much more involved in comparison to the
one-dimensional setting in Section \ref{boundary}.

\subsection{SPDEs with spatially homogeneous noise}
\label{spatial}

Let us consider the following stochastic parabolic Cauchy problem in
$\mathbb{R}^d$: \begin{equation} \frac{\partial u}{\partial
t}(t,x)-\Delta u(t,x)=b(u(t,x)) + \sigma \dot W(t,x),\quad (t,x)\in
[0,T]\times \mathbb{R}^d,\label{24} \end{equation} where $T>0$,
$\Delta$ denotes the Laplacian operator in $\mathbb{R}^d$,
$b:\mathbb{R}\rightarrow \mathbb{R}$ is a $\mathcal{C}^1$ function
with bounded derivative, and suppose that we are given an initial
condition of the form
$$u(0,x)=u_0(x),\; x\in \mathbb{R}^d,$$
with $u_0:\mathbb{R}^d\rightarrow \mathbb{R}$ measurable and
bounded.

The random perturbation $\dot W$ (formally) stands for a Gaussian
noise which is white in time and with some spatially homogeneous
correlation. More precisely, on a complete probability space
$(\Omega, \mathcal{F},P)$, we consider a family of mean zero
Gaussian random variables $W = \{W(\varphi),\, \varphi \in
\mathcal{C}_0^\infty (\mathbb{R}^{d+1})\}$, where
$\mathcal{C}_0^\infty (\mathbb{R}^{d+1})$ denotes the space of
infinitely differentiable functions with compact support, with
covariance functional \begin{equation} E(W(\varphi) W(\psi)) =
\int_0^\infty \int_{\mathbb{R}^d} \left(
\varphi(t)*{\psi_s}(t)\right) (x) \, \Lambda(dx) dt, \; \varphi,
\psi\in \mathcal{C}^\infty_0(\mathbb{R}^{d+1}), \label{25}
\end{equation} where $\psi_s(t,x):=\psi(t,-x)$ and $\Lambda$ is a
non-negative and non-negative definite tempered measure. By
\cite[Chapter VII, Th\'eor\`eme XVII]{sch}, $\Lambda$ is symmetric
and there exists a non-negative tempered measure $\mu$ whose Fourier
transform is $\Lambda$. That is, by definition of the Fourier
transform on the space $\mathcal{S}'(\mathbb{R}^d)$ of tempered
distributions, for all $\phi$ belonging to the space
$\mathcal{S}(\mathbb{R}^d)$ of rapidly decreasing
$\mathcal{C}^{\infty}$ functions,
$$\int_{\mathbb{R}^d} \phi(x) \Lambda(dx)=\int_{\mathbb{R}^d} \mathcal{F} \phi(\xi) \mu(d\xi),$$
and there is an integer $m\geq 1$ such that
$$\int_{\mathbb{R}^d} (1+|\xi|^2)^{-m}\mu(d\xi) <\infty.$$
The measure $\mu$ is usually called the spectral measure of the
noise $W$. In particular, we have
$$  E(W(\varphi) W(\psi)) = \int_0^\infty  \int_{\mathbb{R}^{d}}
\mathcal{F} \varphi(t)(\xi) \overline{\mathcal{F} \psi(t)(\xi)} \,
\mu(d\xi) dt,\; \varphi, \psi\in \mathcal{C}_0^\infty
(\mathbb{R}^{d+1}).$$ Notice that we have used the symbol ``$*$''
for the standard convolution in $\mathbb{R}^d$.
\begin{example}
Usual examples of spatial correlations are given by
$\Lambda(dx)=f(x)dx$, where $f$ is a non-negative continuous
function on $\mathbb{R}^{d}\setminus \{0\}$ which is integrable in a
neighborhood of $0$; for instance, one can take $f$ to be a Riesz
kernel $f_\epsilon(x)=|x|^{-\epsilon}$, for $0<\epsilon<d$. The
space-time white noise would correspond to consider $f$ equal to the
Dirac delta on the origin. In this latter case, the spectral measure
$\mu$ is the Lebesgue measure on $\mathbb{R}^d$.
\end{example}

We denote by $\mathcal{H}$ the completion of the Schwartz space
$\mathcal{S} (\mathbb{R}^d)$ endowed with the semi-inner product
$$\langle \phi_1,\phi_2 \rangle_{\mathcal{H}} =\int_{\mathbb{R}^d} \left(\phi_1 * (\phi_2)_s \right) (x) \, \Lambda(dx)  =\int_{\mathbb{R}^d}
 \mathcal{F} \phi_1(\xi) \overline{\mathcal{F} \phi_2(\xi)} \, \mu(d\xi),$$
$\phi_1,\phi_2\in \mathcal{S}(\mathbb{R}^d)$, and associated
semi-norm $\| \cdot \|_{\mathcal{H}} $. Notice that $\mathcal{H}$ is
a Hilbert space that may contain distributions (see \cite[Example
6]{Da}). Set $\mathcal{H}_T:=L^2([0,T];\mathcal{H})$.

Then, it turns out that the Gaussian noise $W$ can be naturally
extended to $\mathcal{H}_T$, so that we obtain a family $\{W(h),\;
h\in \mathcal{H}_T\}$ of centered Gaussian random variables such
that
$$E(W(h_1)W(h_2))= \langle h_1,h_2\rangle_{\mathcal{H}_T}= \int_0^T \langle h_1(t),h_2(t)\rangle_\mathcal{H}\, dt,\;
h_1,h_2\in \mathcal{H}_T.$$ The details of that extension can be
found, for instance, in \cite{MS}, p. 805, or \cite{QS}, p.3. Then,
if we set
$$W_t(g):=W(\textbf{1}_{[0,t]} g),\; t\in [0,T], g\in \mathcal{H},$$
we obtain a cylindrical Wiener process $\{W_t(g),\; t\in [0,T],g\in
\mathcal{H}\}$ on the Hilbert space $\mathcal{H}$. That is, for any
$t\in [0,T]$ and $g\in \mathcal{H}$, $W_t(g)$ is a mean zero
Gaussian random variable and
$$E(W_t(g_1)W_s(g_2))=(t\wedge s)\langle g_1,g_2\rangle_\mathcal{H},\; s,t\in [0,T],\; g_1,g_2\in \mathcal{H}.$$
We denote by $\mathcal{F}_t$ the $\sigma$-algebra generated by
$\{W_s(g),\; s\in [0,t],g\in \mathcal{H}\}$ and the $P$-null sets.
As it has been explained in \cite[Section 3]{nq}, one can construct
(real-valued) stochastic integrals of predictable processes in
$L^2(\Omega\times [0,T];\mathcal{H})$ with respect to the
cylindrical Wiener process $W$; notice that we are making an abuse
of notation since $W$ denoted the Gaussian family given at the
beginning. The resulting stochastic integral turns out to extend
Walsh's integration theory \cite{walsh} and it is equivalent, in
some particular situations, to Dalang's stochastic integral set up
in \cite{Da}. The stochastic integrals appearing throughout this
section will be understood as integrals with respect to the
cylindrical Wiener
process $W$.\\

We are now in position to rigorously define the mild solution of
Equation (\ref{24}): a square integrable $\{\mathcal{F}_t\}$-adapted
random field $\{u(t,x),\; (t,x)\in [0,T]\times \mathbb{R}^d\}$
solves Equation (\ref{24}) if, for any $(t,x)\in [0,T]\times
\mathbb{R}^d$,
\begin{align}
u(t,x) = & \int_{\mathbb{R}^d} G_t(x-y) u_0(y)\, dy + \int_0^t
\int_{\mathbb{R}^d} G_{t-s}(x-y)
b(u(s,y)\, dy ds\nonumber \\
&\quad \quad + \sigma \int_0^t \int_{\mathbb{R}^d} G_{t-s}(x-y)
W(ds,dy). \label{26}
\end{align}
Here, the function $G_t(x)$, $t>0$, $x\in \mathbb{R}^d$, denotes the
fundamental solution associated to the heat equation in
$\mathbb{R}^d$, that is the centered Gaussian kernel of variance
$2t$:
$$G_t(x)=\frac{1}{\sqrt{4\pi t}}e^{-\frac{x^2}{4t}}.$$
General existence and uniqueness results for Equation (\ref{26}) may
be found in \cite{Da}. More precisely, it turns out
that sufficient conditions for existence and uniqueness of solution
to Equation (\ref{26}) are the following: $b$ is Lipschitz, $u_0$ is
measurable and bounded, and the noise's spatial correlation is
related with the fundamental solution $G$ through the condition
\begin{equation} \int_0^T \int_{\mathbb{R}^d} |\mathcal{F} G_t(\xi)|^2 \; \mu(d\xi)
dt<+\infty. \label{27} \end{equation} As it has been shown in
\cite[Example 8]{Da}, condition (\ref{27}) holds if and only if
\begin{equation} \int_{\mathbb{R}^d} \frac{1}{1+|\xi|^2}\mu(d\xi)<+\infty,
\label{28} \end{equation} and this integrability condition will
be assumed to be satisfied in the remainder of this section.\\

Let us turn now to the question whether the probability law of the
solution at any point has a density. The Gaussian setting in which
we apply the Malliavin calculus machinery is determined by the
Gaussian family $\{W(h),\; h\in \mathcal{H}_T\}$ given before.

Using this framework, it is a consequence of \cite[Theorem 5.2]{nq}
that, if the drift coefficient $b\in \mathcal{C}^1$ has a bounded
and Lipschitz continuous derivative and condition (\ref{28}) is
fulfilled, then the solution $u(t,x)$ to Equation (\ref{26}), at any
point $(t,x)\in (0,T]\times \mathbb{R}^d$, is differentiable in the
Malliavin sense, that is $u(t,x)\in \mathbb{D}^{1,2}$, and its law
has a density with respect to the Lebesgue measure. At this point,
we should mention that in \cite{nq} all the results are proved in
the case where the noise's correlation is given by a function $f$,
that is $\Lambda(dx)=f(x)dx$. However, the extension of those
results to a general tempered measure $\Lambda$ is straightforward.
See also the works \cite{mms}, \cite{QS}, \cite{qs2} for related results with
slightly stronger conditions on the spectral measure.

It will be of much importance for us the equation satisfied be the
Malliavin derivative of $u(t,x)$. Indeed, see either \cite{mms},
\cite{QS} or \cite{nq}, the Malliavin derivative $Du(t,x)$ takes
values in the Hilbert space $\mathcal{H}_T$ and satisfies the
following linear parabolic equation: \begin{equation} Du(t,x)=
\sigma G_{t-\cdot}(x-\star) + \int_0^t \int_{\mathbb{R}^d}
G_{t-s}(x-y) b'(u(s,y)) Du(s,y) \; dy ds, \label{29} \end{equation}
where ``$\star$'' stands for the $\mathcal{H}$-variable. Moreover,
one proves that \begin{equation} \sup_{t\in [0,T]} \sup_{x\in
\mathbb{R}^d} E(\|Du(t,x)\|_{\mathcal{H}_T}^2)<+\infty.
\label{46}\end{equation} Equation (\ref{29}) may be interpreted in
the following sense: for any $r\in [0,t)$, $D_ru(t,x)$ satisfies the
equation in $\mathcal{H}$
\begin{equation} D_ru(t,x)= \sigma G_{t-r}(x-\star) + \int_0^t
\int_{\mathbb{R}^d} G_{t-s}(x-y) b'(u(s,y)) D_ru(s,y) \; dy ds.
\label{30}
\end{equation}
The integral on the right-hand side of (\ref{30}) is understood
as a $\mathcal{H}$-valued pathwise integral. Before going on with
our analysis, let us briefly describe how this integral is
rigorously defined. Let $\{e_j,\; j\geq 1\}$ be a complete orthonormal
system of $\mathcal{H}$. Then, using the properties of $G$, the
boundedness of $b'$ and (\ref{46}), the $\mathcal{H}$-valued
integral
$$\mathcal{I}_t=\int_0^t  \int_{\mathbb{R}^{d}} G_{t-s}(x-y) b'(u(s,y)) D_ru(s,y)\; dy ds$$
can be defined through its components $$\left\{ \int_0^t
\int_{\mathbb{R}^{d}} G_{t-s}(x-y) b'(u(s,y)) \langle
D_ru(s,y),e_j\rangle_{\mathcal{H}}\; dy ds, \;j\geq 1\right\}$$ with
respect to the basis $\{e_j,\; j\geq 1\}$, and those latter integrals
takes values in $\mathbb{R}$. Moreover, one can obtain an upper
bound for the square moment of $\mathcal{I}_t$ (for the general setting see, for instance,
\cite{quer}, p. 24):
\begin{equation}
E(|\mathcal{I}_t|^2)\leq C \int_0^t \int_{\mathbb{R}^d}
|G_{t-s}(x-y)|^2 |b'(u(s,y))|^2 E(\|D_r u(s,y)\|_{\mathcal{H}}^2) \;
dy ds. \label{bound-int}
\end{equation}

Let us go back now to Equation (\ref{30}). Since
$G_{t-r}(x-\star):\mathbb{R}^d\rightarrow \mathbb{R}$ defines a
function (indeed in $\mathcal{S}(\mathbb{R}^d)$), this implies that
either $D_ru(t,x)$ and the integral in the right-hand side of
(\ref{30}) define elements of $\mathcal{H}$ which are functions in
$z$. Therefore, for any fixed $(t,x)\in (0,T]\times \mathbb{R}^d$,
we can state that the Malliavin derivative satisfies, for all
$(r,z)\in [0,t)\times \mathbb{R}^d$:
\begin{equation} D_{r,z}u(t,x)=  \sigma G_{t-r}(x-z) + \int_0^t \int_{\mathbb{R}^d}
G_{t-s}(x-y) b'(u(s,y)) D_{r,z}u(s,y) \; dy ds. \label{31}
\end{equation} A crucial consequence of this fact is that, as in the
case of the one-dimensional stochastic heat equation with boundary
conditions, the Malliavin derivative $D_{r,z}u(t,x)$ is
non-negative, for all
$(r,z)\in [0,t)\times \mathbb{R}^d$, a.s.\\

We will need a slightly stronger condition on the spectral measure
$\mu$ than (\ref{28}). Namely, consider the following hypothesis:\\

\noindent \textbf{Hypothesis H$_\eta$:} There exists $\eta\in (0,1)$
such that
\begin{equation}
\int_{\mathbb{R}^d} \frac{1}{(1+|\xi|^2)^\eta} \mu(d\xi)<+\infty.
\label{69}
\end{equation}

\vspace{0.3cm}

Then, one can prove the following estimates (see \cite[Lemma
3.1]{mms}):

\begin{lemma}\label{lemma4}
Assume that the spectral measure $\mu$ satisfies
$\int_{\mathbb{R}^d} \frac{1}{1+|\xi|^2}\mu(d\xi)<+\infty$.
\begin{enumerate}
  \item Let $T>0$. Then, there exists a constant $k_1>0$ such that, for all $t\in [0,T]$,
$$k_1 t \leq \int_0^t \int_{\mathbb{R}^d} |\mathcal{F} G_s(\xi)|^2 \;
\mu(d\xi) ds. $$ The constant $k_1$ depends on $T$ and, indeed, it
converges to zero as $T$ tends to infinity.
  \item Suppose that Hypothesis H$_\eta$ holds. Then, there exists a
  constant $k_2>0$ such that, for all $t\geq 0$,
  \begin{equation} \int_0^t \int_{\mathbb{R}^d}
|\mathcal{F} G_s(\xi)|^2 \; \mu(d\xi) ds\leq k_2 t^\beta,
\label{32}\end{equation} for all $\beta\in (0,1-\eta]$.
\end{enumerate}
\end{lemma}

\begin{remark}\label{remark3}
It is worth mentioning that the integrability condition (\ref{27})
was sufficient for us  to prove the existence of density for the
solution $u(t,x)$, at any point $(t,x)\in (0,T]\times \mathbb{R}^d$
(see \cite[Theorem 5.2]{nq}). However, as it will be made clearer in
Section \ref{main-spatial}, we will really need upper bounds of the
form (\ref{32}) in order to obtain Gaussian estimates for the
density of $u(t,x)$.
\end{remark}

\subsection{Gaussian estimates for the density of the solution}
\label{main-spatial}

Let us consider $T>0$ and $\{u(t,x),\; (t,x)\in [0,T]\times
\mathbb{R}^d\}$ the unique mild solution to Equation (\ref{26}).
This section is devoted to proof the following result:

\begin{theorem}
Fix $t\in (0,T]$ and $x\in \mathbb{R}^d$. Suppose that Hypothesis
H$_\eta$ is satisfied for some $\eta\in (0,\frac{3}{4})$ and that
the coefficient $b$ is of class $\mathcal{C}^1$ and has a bounded
Lipschitz continuous derivative. Then, the random variable $u(t,x)$
has a density $p$ with respect to Lebesgue measure which satisfies
the following: for almost every $z\in \mathbb{R}$, $$\frac{E\left|
u(t,x)-m\right|}{C_2 t^{1-\eta}} \exp\left\{ -\frac{(z-m)^2}{ C_1 t
}\right\}\leq p(z) \leq \frac{E\left| u(t,x)-m\right|}{ C_1 t}
\exp\left\{ -\frac{(z-m)^2}{ C_2 t^{1-\eta}}\right\}, $$ where
$m=E(u(t,x))$ and $C_1$, $C_2$ are positive constants   depending on
$\sigma$, $\|b'\|_\infty$, $\eta$ and $T$. \label{bounds}
\end{theorem}


Theorem \ref{bounds} will be a consequence of \cite[Theorem 3.1]{nv}
and the following proposition. As we have done in Section
\ref{main}, we use the notation $F=u(t,x)-E(u(t,x))$ and we remind
that we will need to find almost sure lower and upper bounds for the
random variable $g(F)$, where
$$g(F) =\int_0^\infty e^{-\theta} E\left[ E'\left( \langle
Du(t,x),\widetilde{Du(t,x)}\rangle_{\mathcal{H}_T}\right)\big|
F\right] d\theta.$$

\begin{proposition}
Fix $T>0$. Assume that Hypothesis H$_\eta$ holds for some $\eta\in
(0,\frac{3}{4})$ and the coefficient $b$ is of class $\mathcal{C}^1$
and has a bounded Lipschitz continuous derivative. There exist
positive constants $C_1, C_2$ such that,
\begin{equation} C_1 t\leq
g(F)\leq C_2 t^{1-\eta},\; a.s., \label{36} \end{equation} for all
$t\in (0,T]$. \label{prop2}
\end{proposition}

In order to prove Proposition \ref{prop2}, we will need the
following technical lemma, which plays the role of Lemma
\ref{lemma1} in our standing setting.

\begin{lemma}
Let $t>0$ and assume that Hypothesis H$_\eta$ holds. Then, there
exists a positive constant $C$  depending on $\sigma$, $\|b'\|_\infty$ and the constant $k_2$    in (\ref{32}), such that, for all $\delta\in (0,1]$,
\begin{equation} \sup_{\stackrel{(1-\delta)t\leq \nu\leq t}{y\in \mathbb{R}^d}} E \left[
\int_{(1-\delta)t}^t \| D_r u(\nu,y)\|_{\mathcal{H}}^2\; dr \Big| F
\right] \leq C (\delta t)^\beta,\; a.s. \label{33} \end{equation}
and
\begin{equation} \sup_{\theta \geq 1}
\sup_{\stackrel{(1-\delta)t\leq \nu\leq t}{y\in \mathbb{R}^d}} E
\left[ E'\left(\int_{(1-\delta)t}^t \|
\widetilde{D_ru(\nu,y)}\|_{\mathcal{H}}^2 \; dr \right) \Big| F
\right] \leq C (\delta t)^\beta, \; a.s., \label{34}
\end{equation} for any $\beta \in (0,1-\eta]$. \label{lemma2}
\end{lemma}
\emph{Proof:} It is very similar to that of Lemma \ref{lemma1}. In
fact, owing to Equation (\ref{30}), we have, for any $(\nu,v)\in
[(1-\delta)t,t]\times \mathbb{R}^d$,
\begin{align}
& \int_{(1-\delta)t}^t \|D_ru(\nu,v)\|_\mathcal{H}^2 \;dr \leq  2
\sigma^2
\int_{(1-\delta)t}^t \|G_{\nu-r}(v-\star)\|_\mathcal{H}^2 \;dr \nonumber \\
&\quad \quad + 2 \|b'\|_\infty^2 (\delta t) \int_{(1-\delta)t}^\nu
\int_{\mathbb{R}^d} |G_{\nu-s}(x-y)|^2 \left(\int_{(1-\delta)t}^t
\|D_ru(s,y)\|_\mathcal{H}^2\; dr \right) \; dy ds, \label{35}
\end{align}
where we have applied Minkowski's and Cauchy-Schwarz's inequalities.
By (\ref{32}) in Lemma \ref{lemma4},
$$\int_{(1-\delta)t}^t \|G_{\nu-r}(v-\star)\|_\mathcal{H}^2 \; dr =
\int_0^{\delta t} \int_{\mathbb{R}^d} |\mathcal{F}
G_r(\xi)|^2\;\mu(d\xi) dr\leq k_2 (\delta t)^\beta,$$ for all
$\beta\in (0,1-\eta]$. Therefore, plugging this bound in (\ref{35})
and taking conditional expectation, we obtain:
\begin{align*}
& E\left[ \int_{(1-\delta)t}^t \|D_ru(\nu,v)\|_\mathcal{H}^2 \;dr
\Big| F\right]
\leq   2\sigma^2 k_2  (\delta t)^\beta \nonumber \\
&\quad \quad +  \frac{1}{2 \sqrt{2\pi}}\, \|b'\|_\infty^2  (\delta
t) \int_{(1-\delta)t}^\nu \left( \sup_{\stackrel{(1-\delta)t \leq
\tau \leq s}{y\in \mathbb{R}^d}} E\left[ \int_{(1-\delta)t}^t
\|D_ru(\tau,y)\|_\mathcal{H}^2\; dr \Big| F\right]\right)
\frac{1}{\sqrt{\nu-s}}\; ds, \; a.s.
\end{align*}
As we have done in the proof of Lemma \ref{lemma1}, we are now in
position to apply Gronwall's lemma \cite[Lemma 15]{Da}. Hence
(\ref{33}) is proved. The estimation (\ref{34}) can be checked using
exactly the same arguments. \qed

\vspace{0.3cm}

\emph{Proof of Proposition \ref{prop2}:} The framework of the proof
is similar to that of Proposition \ref{prop} in Section \ref{main}.
However, computations here will be slightly more involved since we
are working in a Hilbert-space-valued setting determined by
$\mathcal{H}_T=L^2([0,T];\mathcal{H})$. Let us first deal with the
lower bound in
(\ref{36}):\\

\emph{Step 1: The lower bound.} Recall that $F=u(t,x)-E(u(t,x))$ and
the random variable $g(F)$ can be written as
\begin{align*}
g(F)& =\int_0^\infty e^{-\theta} E\left[ E'\left( \langle
Du(t,x),\widetilde{Du(t,x)}\rangle_{\mathcal{H}_T}\right)\big|
F\right]
d\theta \\
& =\int_0^\infty e^{-\theta} E\left[ E'\left( \int_0^t \langle
D_ru(t,x),\widetilde{D_ru(t,x)}\rangle_\mathcal{H} \; dr
\right)\Big| F\right] d\theta,
\end{align*}
where $\widetilde{Du(t,x)}$ denotes the shifted random variable
$(Du(t,x))(e^{-\theta}\omega +\sqrt{1-e^{-2\theta}}\omega')$.

According to Equation (\ref{29}), for any $\delta \in (0,1]$ we have
the decomposition
 \begin{equation} g(F)\geq  \sigma^2 B_0(t,x;\delta)-
   \left|B_1(t,x;\delta) + B_2(t,x;\delta) + B_3(t,x;\delta )\right|,
\label{37} \end{equation} where

$$B_0(t,x;\delta)= \int_{(1-\delta)t}^t \|G_{t-r}(x-\star)\|_\mathcal{H}^2 \; dr,$$
\begin{align*}
& B_1(t,x;\delta) \\
& \; = \sigma E\left[ \int_{(1-\delta)t}^t \left\langle
G_{t-r}(x-\star),\int_{(1-\delta)t}^t \int_{\mathbb{R}^d}
G_{t-s}(x-y) b'(u(s,y)) D_ru(s,y)\; dy ds\right\rangle_\mathcal{H}
dr \Big| F\right],
\end{align*}
\begin{align*}
B_2(t,x;\delta) = &\sigma \int_0^\infty e^{-\theta} E\left[ E'\left(
\int_{(1-\delta)t}^t
\left\langle G_{t-r}(x-\star) \right. \right. \right.\\
& \quad ,\left.\left.\left. \int_{(1-\delta)t}^t \int_{\mathbb{R}^d}
G_{t-s}(x-y) b'(\widetilde{u(s,y)}) \widetilde{D_ru(s,y)}\; dy
ds\right\rangle_\mathcal{H} dr \right) \Big| F\right] d\theta,
\end{align*}
\begin{align*}
  B_3(t,x;\delta)= &  \int_0^\infty e^{-\theta} E\left[ E'\left( \int_{(1-\delta)t}^t \left\langle
  \int_{(1-\delta)t}^t \int_{\mathbb{R}^d} G_{t-s}(x-y) b'(u(s,y)) D_ru(s,y)\; dyds
  \right.\right.\right. \\
  & \quad \quad  \left. \left.\left., \int_{(1-\delta)t}^t \int_{\mathbb{R}^d} G_{t-s}(x-y) b'(\widetilde{u(s,y)}) \widetilde{D_ru(s,y)}\;
  dyds \right\rangle_\mathcal{H} dr \right) \Bigg| F\right] d\theta.
\end{align*}
By part 1 in Lemma \ref{lemma4}, notice first that \begin{equation}
B_0(t,x;\delta)=\int_0^{\delta t} \int_{\mathbb{R}^d} |\mathcal{F}
G_s(\xi)|^2\;\mu(d\xi) ds\geq  k_1 \; \delta t. \label{38}
\end{equation}

Concerning the second term $B_1(t,x;\delta)$, we can apply
Cauchy-Schwarz and Mikowski's inequalities, so that we obtain
\begin{align*}
|B_1(t,x;\delta)|\leq & C \; (\delta t)^{\frac{1}{2}}
\left( \int_{(1-\delta)t}^t \|G_{t-r}(x-\star)\|_\mathcal{H}^2 dr\right)^{\frac{1}{2}} \nonumber \\
& \; \times \left( E\left[ \int_{(1-\delta)t}^t \left\|
\int_{(1-\delta)t}^t \int_{\mathbb{R}^d} G_{t-s}(x-y) b'(u(s,y))
D_ru(s,y) \; dy ds \right\|_\mathcal{H}^2 dr \Bigg| F\right]
\right)^{\frac{1}{2}}
\nonumber \\
& \leq C \; (\delta t)^{\frac{1}{2}}
\left( \int_0^{\delta t}  \int_{\mathbb{R}^d}  |\mathcal{F} G_r (\xi)|^2 \; \mu(d\xi) dr\right)^{\frac{1}{2}} \nonumber \\
& \quad \times \left( \int_{(1-\delta)t}^t \int_{\mathbb{R}^d}
|G_{t-s}(x-y)|^2 E\left[ \int_{(1-\delta)t}^t
\|D_ru(s,y)\|_\mathcal{H}^2 dr \Big| F\right] dy
ds\right)^{\frac{1}{2}}.
\end{align*} Thus, by (\ref{32}), Lemma \ref{lemma2} and the fact
that \begin{equation} \int_{(1-\delta)t}^t \int_{\mathbb{R}^d}
|G_{t-s}(x-y)|^2 \; dy ds \leq C (\delta t)^{\frac{1}{2}},
\label{43}\end{equation} we have:
 \begin{equation} |B_1(t,x;\delta)|\leq C (\delta t)^{\beta+\frac{3}{4}}, \label{39} \end{equation}
for all $\beta\in (0,1-\eta]$. The term $|B_2(t,x;\delta)|$ can be
treated in the same way as we have just done for
$|B_1(t,x;\delta)|$. Namely, one proves that
\begin{align*}
& |B_2(t,x;\delta)|\leq  C \; (\delta t)^{\frac{1}{2}}
\left( \int_{(1-\delta)t}^t \|G_{t-r}(x-\star)\|_\mathcal{H}^2 dr\right)^{\frac{1}{2}} \nonumber \\
& \; \times \int_0^\infty e^{-\theta} \left( E\left[ E'\left(
\int_{(1-\delta)t}^t \left\| \int_{(1-\delta)t}^t
\int_{\mathbb{R}^d} G_{t-s}(x-y) b'(u(s,y)) D_ru(s,y) \; dy ds
\right\|_\mathcal{H}^2 dr \right) \Bigg| F\right]
\right)^{\frac{1}{2}} d\theta
\nonumber \\
& \leq C \; (\delta t)^{\frac{1}{2}}
\left( \int_0^{\delta t}  \int_{\mathbb{R}^d}  |\mathcal{F} G_r (\xi)|^2 \; \mu(d\xi) dr\right)^{\frac{1}{2}} \nonumber \\
& \quad \times \int_0^\infty e^{-\theta} \left( \int_{(1-\delta)t}^t
\int_{\mathbb{R}^d} |G_{t-s}(x-y)|^2 E\left[ E'\left(
\int_{(1-\delta)t}^t \|D_ru(s,y)\|_\mathcal{H}^2 dr \right) \Big|
F\right] dy ds\right)^{\frac{1}{2}} d\theta.
\end{align*}
Taking into account (\ref{32}), Lemma \ref{lemma2} and (\ref{43}),
we also get \begin{equation} |B_2(t,x;\delta)|\leq C (\delta
t)^{\beta+\frac{3}{4}}, \label{40} \end{equation} for all $\beta\in
(0,1-\eta])$.

Eventually, in order to deal with the term $|B_3(t,x;\delta)|$, we
mainly apply Cauchy-Schwarz's inequality with respect to
$E\left[E'\left( \int_{(1-\delta)t}^t \| \bullet \|_\mathcal{H}^2 dr
\right)\Big| F\right]$, so that we will take advantage of the
computations which we have performed so far. More precisely, we have

\begin{align*}
  & B_3(t,x;\delta) \\
  & \;\leq  \left( E\left[ \int_{(1-\delta)t}^t \left\|
\int_{(1-\delta)t}^t \int_{\mathbb{R}^d} G_{t-s}(x-y) b'(u(s,y))
D_ru(s,y)
\; dy ds \right\|_\mathcal{H}^2 dr \Bigg| F\right] \right)^{\frac{1}{2}} \\
&\; \times \int_0^\infty e^{-\theta} \left( E\left[ E'\left(
\int_{(1-\delta)t}^t \left\| \int_{(1-\delta)t}^t
\int_{\mathbb{R}^d} G_{t-s}(x-y) b'(u(s,y)) D_ru(s,y) \; dy ds
\right\|_\mathcal{H}^2 dr \right) \Bigg| F\right]
\right)^{\frac{1}{2}} d\theta.
\end{align*}
The two terms in the right-hand side of the above inequality already
appeared in the analysis of $B_1(t,x;\delta)$ and $B_2(t,x;\delta)$,
respectively, and each of them may be bounded, up to some constant,
by $(\delta t)^{\frac{\beta}{2}+\frac{3}{4}}$. Therefore
\begin{equation} |B_3(t,x;\delta)|\leq C (\delta
t)^{\beta+\frac{3}{2}}, \label{41} \end{equation} for any
$\beta\in (0,1-\eta]$. Estimations (\ref{37}), (\ref{38}) and (\ref{39})-(\ref{41}) yield
$$
g(F)\geq  \sigma^2 k_1 \delta t  - c_1 \left(
(\delta t)^{\beta+\frac{3}{4}}+ (\delta
t)^{\beta+\frac{3}{2}}\right),
$$
where $c_1$ depends on $\sigma$, $\|b'\|_\infty$ and $k_2$.
Hence, if $\delta< 1\wedge
\frac{1}{T}$ and $\beta \in (\frac{1}{4},1-\eta]$,
$$g(F)\geq  t \left(    \sigma^2 k_1 \delta - 2 c_1 \delta^{\beta+\frac{3}{4}}  t^{\beta-\frac{1}{4}} \right)\geq
t \left(    \sigma^2 k_1 \delta - 2 c_1 \delta^{\beta+\frac{3}{4}}
T^{\beta-\frac{1}{4}} \right).$$ Observe that the quantity
$$C_1:= \sigma^2 k_1 \delta - 2 c_1 \delta^{\beta+\frac{3}{4}}
T^{\beta-\frac{1}{4}}$$ defines a positive constant whenever
$\delta\in (0,\delta_0)$, with
$$\delta_0=1\wedge \frac{1}{T}\wedge
\frac{1}{T}\left(\frac{ \sigma^2 k_1}{2c_1}\right)^{\frac{1}{\beta-\frac{1}{4}}}.$$
Therefore, we obtain the desired lower bound in (\ref{36}).\\

\emph{Step 2: The upper bound.} The upper in (\ref{36}) is an almost
immediate consequence of the computations in the Step 1 and the
decomposition \begin{equation} g(F)\leq \sum_{i=0}^3
|B_i(t,x;1)|.\label{42}\end{equation} Indeed, we have already found
upper bounds for the last three terms on the right-hand side of
(\ref{42}). On the other hand, observe that (\ref{32}) yields
$$B_0(t,x)\leq C t^\beta,$$
for all $\beta\in (0,1-\eta]$. This bound, together with
(\ref{39})-(\ref{41}) in the case $\delta=1$, implies
$$g(F)\leq C \left( t^\beta + t^{\beta + \frac{3}{4}} + t^{\beta
+\frac{3}{2}}\right)\leq C_2 t^\beta,$$ for all $\beta\in
(0,1-\eta]$, where the constant $C_2$ depends on $T$. Therefore we
conclude the proof. \qed

\vspace{0.3cm}

The proof of Theorem \ref{bounds} can be finished as in the case of Theorem \ref{teo1}.
 \qed


\section{Gaussian bounds for the stochastic wave equation at small time}
\label{wave}

The main objective here is to extend the results in Section \ref{initial} to a stochastic wave equation in space dimension $d\leq 3$ and controlled
by the spatially homogeneous Gaussian noise considered there. The intrinsic properties of the differential operator driving the equation will not allow us to obtain optimal results for all time parameter $T>0$, even if we assume that the noise's space correlation satisfies Hypothesis H$_\eta$, an assumption which is slightly stronger than the one needed to have existence of density for the corresponding mild solution (see \cite{nq}).

\subsection{The stochastic wave equation in dimension $d=1,2,3$}
\label{setting-wave}

We consider here the same setting as in Section \ref{initial} but for the stochastic wave equation in spatial dimension $d\leq 3$:
\begin{equation}
\frac{\partial^2 u}{\partial
t^2}(t,x)-\Delta u(t,x)=b(u(t,x)) + \sigma \dot W(t,x),\quad (t,x)\in
[0,T]\times \mathbb{R}^d,
\label{50}
\end{equation}
where $T>0$,
$b:\mathbb{R}\rightarrow \mathbb{R}$ is a $\mathcal{C}^1$ function
with bounded derivative, and suppose that we are given initial
conditions of the form
$$u(0,x)=u_0(x),\; \frac{\partial u}{\partial t}(0,x)=v_0(x),\;   x\in \mathbb{R}^d,$$
with $u_0,v_0:\mathbb{R}^d\rightarrow \mathbb{R}$ measurable and
bounded functions such that $u_0$ is of class
$\mathcal{C}^1(\mathbb{R}^d)$ and has a bounded derivative $\nabla
u_0$. The random perturbation $\dot W$ corresponds to the spatially
homogeneous Gaussian noise described in the previous
Section \ref{spatial}. We recall that $\mu$ denotes the corresponding spectral measure and $\{\mathcal{F}_t,\;t\geq 0\}$ the filtration defined by the cylindrical Wiener process associated to the noise $W$.\\

The mild solution of Equation (\ref{50}) is given by a $\{\mathcal{F}_t\}$-adapted process $\{u(t,x),\; (t,x)\in [0,T]\times \mathbb{R}^d\}$ such that,
for all $(t,x)\in (0,T]\times \mathbb{R}^d$,
\begin{align}
u(t,x) = & \int_{\mathbb{R}^d} v_0(x-y) \Gamma^d_t(dy) + \frac{\partial}{\partial t} \left(\int_{\mathbb{R}^d} u_0(x-y) \Gamma^d_t(dy)\right) \nonumber \\
&\; + \int_0^t
\int_{\mathbb{R}^d}
b(u(s,x-y)\, \Gamma^d_{t-s}(dy) ds + \sigma \int_0^t \int_{\mathbb{R}^d} \Gamma^d_{t-s}(x-y)
W(ds,dy),
\label{51}
\end{align}
where $\Gamma_t^d$, $t> 0$, denotes the fundamental solution of the wave equation in dimension $d=1,2,3$:
\begin{align*}
\Gamma_t^1(x) &=\frac{1}{2} {\bf 1}_{\{|x|<t\}},\\
\Gamma_t^2(x) & = \frac{1}{2\pi} (t^2-|x|^2)_+^{-1/2},\\
\Gamma_t^3 &= \frac{1}{4\pi t}\sigma_t,
\end{align*}
The element $\sigma_t$ stands for the surface measure on the
three-dimensional sphere of radius $t$. In particular, for each $t$,
$\Gamma^d_t$ has compact support and, in the case $d=3$,
$\Gamma^3_t$ is no more a function but measure on $\mathbb{R}^3$. It
is important to remark that only in these cases, $\Gamma^d$,
$d=1,2,3$, defines a non-negative measure. Existence and uniqueness
of mild solution to Equation (\ref{51}) is a consequence of the
results in \cite{dq}, whenever the space correlation satisfies
\begin{equation}
\int_{\mathbb{R}^d}\frac{1}{1+|\xi|^2}\mu(d\xi)<+\infty. \label{54}
\end{equation} We also point out that the stochastic integral in the
right-hand side of (\ref{51}) is a well-defined integral of a
deterministic element in $\mathcal{H}_T$ with respect to the
cylindrical Wiener process associated to the noise (see Lemma 3.2
and Example 4.2 in \cite{nq}).

For all dimensions $d\geq 1$, we have a unified expression for the Fourier transform of $\Gamma^d_t$:
$$\mathcal{F} \Gamma^d_t (\xi)= \frac{\sin (2\pi t |\xi|)}{2\pi |\xi|}.$$
Using this fact and assuming that (\ref{54}) holds,  one proves the
following lemma (see \cite{lev}, Lemmas 5.4.1 and 5.4.3):
\begin{lemma}\label{lemma5}
For any $t\geq 0$ it holds that \begin{equation} c_{1}(t\land
t^{3})\frac{1}{1+|\xi|^{2}} \leq \int_{0}^{t} |\mathcal{F}
\Gamma^d_s(\xi)|^{2} ds\leq c_{2} (t+t^{3}) \frac{1}{1+|\xi|^{2}},
\label{70} \end{equation} with some positive constants $c_1,c_2>0$.
\end{lemma}
Thus, if we assume that $ \int_{\mathbb{R}^d}
\frac{\mu(d\xi)}{1+|\xi|^{2}} <\infty$, (\ref{70}) yields
$$
d_{1}(t\land t^{3})\leq \int_{0}^{t} \int_{\mathbb{R}^d}
|\mathcal{F} \Gamma^{d}_s(\xi)|^{2} \mu(d\xi) ds \leq d_{2}
(t+t^{3}),$$ for all $t\geq 0$, with some positive constants $d_{1},
d_2$. In particular, for $t\in [0,1)$ we have
\begin{equation} d_1 t^3 \leq \int_0^t \int_{\mathbb{R}^d} |\mathcal{F}
\Gamma_s^d(\xi)|^2 \; \mu(d\xi) ds\leq d_2 t.
\label{71}\end{equation}
 However, under
Hypothesis H$_\eta$ (see (\ref{69})) one can get a slightly sharper
upper estimation (see \cite[Lemma 3]{QS}):
\begin{lemma}\label{lemma6}
Let $T>0$ and assume that Hypothesis H$_\eta$ holds. Then
\begin{equation} \int_0^t \int_{\mathbb{R}^d} |\mathcal{F}
\Gamma_s^d(\xi)|^2 \; \mu(d\xi) ds\leq d_3 t^{3-2\eta}, \label{52}
\end{equation}
for all $t\in [0,T]$.
\end{lemma}
Eventually, if $d=1,2,3$, explicit computations yield that, for any
$t\geq 0$, \begin{equation} \int_0^t \int_{\mathbb{R}^d} \;
\Gamma^d_s(dy) ds\leq C t^2, \label{53} \end{equation}
where $C$ is a positive constant that only depends on $d$.\\

If $b$ belongs to $\mathcal{C}^1$ and it has a Lipschitz continuous
bounded derivative, then the solution $u(t,x)$, at any $(t,x)\in
(0,T]\times \mathbb{R}^d$, belongs to $\mathbb{D}^{1,2}$ and its
Malliavin derivative, as a random variable taking values in
$\mathcal{H}_T=L^2([0,T];\mathcal{H})$, satisfies
\begin{equation}
Du(t,x)=
\sigma \Gamma^d_{t-\cdot}(x-\star) + \int_0^t \int_{\mathbb{R}^d}
 b'(u(s,x-y)) Du(s,x-y) \; \Gamma_{t-s}^d(dy) ds,
\label{55}
\end{equation}
where ``$\star$'' stands for the $\mathcal{H}$-variable (see
\cite[Proposition 5.1]{nq}). This linear equation is understood in
$L^2(\Omega\times [0,T];\mathcal{H})$ and let us remark that $\sigma
\Gamma^d_{t-\cdot}(x-\star)$ is a well-defined element in
$\mathcal{H}_T$ (see \cite[Lemma 3.2]{nq}). Moreover, under the
standing hypothesis, the random variable $u(t,x)$, for $(t,x)\in
(0,T]\times \mathbb{R}^d$, has a density with respect to the
Lebesgue measure (see \cite[Theorem 5.2]{nq}). Of course, the
Gaussian setting here is the same as the one that has been
considered in Section \ref{initial}.

\subsection{Gaussian estimates of  the density at small time}
\label{bounds-wave}

For $T>0$, consider the unique mild solution $\{u(t,x),\; (t,x)\in
[0,T]\times \mathbb{R}^d\}$ to Equation (\ref{50}). In this section,
we will prove that the density $p:\mathbb{R}^d\rightarrow
\mathbb{R}$ of $u(t,x)$ has lower and upper Gaussian bounds whenever
$T$ is \textit{small}, where this essentially means that $T<1$. The
main result is the following:

\begin{theorem}
Suppose that Hypothesis H$_\eta$ is satisfied and that the
coefficient $b$ is of class $\mathcal{C}^1$ and has a bounded
Lipschitz continuous derivative. Then, there exists $T_0\in (0,1]$
such that the following statement is satisfied: for any $T\in
(0,T_0)$ and $(t,x)\in (0,T]\times \mathbb{R}^d$, the random
variable $u(t,x)$ has a density $p$ with respect to Lebesgue measure
such that, for almost every $z\in \mathbb{R}$,
$$\frac{E\left| u(t,x)-m\right|}{ C_2 t^{3-2\eta}} \exp\left\{ -\frac{(z-m)^2}{C_1 t^3 }
\right\}\leq p(z) \leq \frac{E\left| u(t,x)-m\right|}{ C_1 t^3}
\exp\left\{ -\frac{(z-m)^2}{C_2 t^{3-2\eta}}\right\}, $$ where
$m=E(u(t,x))$ and $C_1, C_2$ are positive constants depending on $\sigma$, $\|b'\|_\infty$, $T_0$ and $\eta$.
 \label{wave-bounds}
\end{theorem}

\begin{remark}
In the case of the stochastic heat equation presented in Section \ref{initial}, we have been able to obtain Gaussian upper and lower bounds  for
any $T>0$, while here we restrict our analysis to \textit{small} $T$. As we will precisely point out in the next Proposition \ref{prop3}, that difference is due to the fact that the Malliavin derivative of the solution to the stochastic wave equation does not need to be a non-negative function.
\end{remark}

The statement of Theorem \ref{wave-bounds} is an immediate
consequence of \cite[Theorem 3.1]{nv} and the following proposition.
For $t>0$ and $x\in \mathbb{R}^d$, set $F=u(t,x)-E(u(t,x))$ and we
remind that we will need to find almost sure lower and upper bounds
for the random variable $g(F)$, where \begin{equation} g(F)
=\int_0^\infty e^{-\theta} E\left[ E'\left( \langle
Du(t,x),\widetilde{Du(t,x)}\rangle_{\mathcal{H}_T}\right)\big|
F\right] d\theta. \label{59} \end{equation}

\begin{proposition}
Assume that Hypothesis H$_\eta$ holds. There exist $T_0\in (0,1]$
and positive constants $C_1,C_2$ such that, for any $T\in (0,T_0)$,
\begin{equation}
C_1 t^3 \leq g(F) \leq C_2 t^{3-2\eta},\; a.s. \label{58}
\end{equation}
for all $t\in [0,T]$. \label{prop3}
\end{proposition}

\textit{Proof:} It follows the same lines as the proof of
Proposition \ref{prop2}, so that we will only point out the main
steps. More precisely, we observe first that the Malliavin
derivative $Du(t,x)$ solves the linear equation (\ref{55}) with a
non-negative initial condition but driven by a hyperbolic operator.
Thus, in comparison with the stochastic heat equation,
$D_{r,z}u(t,x)$ does not need to be non-negative as a function of
$(r,z)$; indeed, in the case $d=3$, even the Malliavin derivative
does not need to be a function. Hence, in order to deal with the
lower bound of $g(F)$ (see (\ref{59})), we will not be able to
restrict the integral with respect to $dr$ on a small time interval
as we have done in the proof of Proposition \ref{prop2}. This is the
reason why we will be forced to consider $T<1$.

In fact, by (\ref{55}), we are only able to consider the
decomposition \begin{equation} g(F)\geq D_0(t)-\left(
|D_1(t)|+|D_2(t)|+|D_3(t)|\right), \label{60} \end{equation} where
$$D_0(t)=\sigma^2 \int_0^t \|\Gamma^d_{t-r}(x-\star)\|_\mathcal{H}^2 \; dr,$$
$$D_1(t)=\sigma E\left[ \int_0^t \left\langle \Gamma^d_{t-r}(x-\star), \int_0^t \int_{\mathbb{R}^d} b'(u(s,x-y)) D_ru(s,x-y) \;
 \Gamma^d_{t-s}(dy) ds \right\rangle_\mathcal{H} dr\Big| F\right].$$
\begin{align*}
D_2(t)= & \int_0^\infty e^{-\theta} \sigma E\left[ E'\left( \int_0^t \left\langle \Gamma^d_{t-r}(x-\star)\right.\right.\right.\\
&\quad , \left.\left.\left.
  \int_0^t \int_{\mathbb{R}^d} b'(\widetilde{u(s,x-y)}) \left(\widetilde{D_ru(s,x-y)}\right) \; \Gamma^d_{t-s}(dy) ds \right\rangle_\mathcal{H} dr\right)\Big| F\right]
d\theta.
\end{align*}

\begin{align*}
D_3(t)= & \int_0^\infty e^{-\theta}  E\left[ E'\left( \int_0^t
\left\langle    \int_0^t \int_{\mathbb{R}^d} b'(u(s,x-y))
D_ru(s,x-y) \;
\Gamma^d_{t-s}(dy) ds    \right.\right.\right.\\
&\quad , \left.\left.\left.
  \int_0^t \int_{\mathbb{R}^d} b'(\widetilde{u(s,x-y)}) \left(\widetilde{D_ru(s,x-y)}\right) \; \Gamma^d_{t-s}(dy) ds \right\rangle_\mathcal{H} dr\right)\Big| F\right]
d\theta.
\end{align*}
By the lower bound in (\ref{71}), we have \begin{equation}
D_0(t)\geq d_1 t^3. \label{61}
\end{equation} Concerning the term $D_1$, we can argue as follows:
\begin{align}
|D_1(t)| & \leq C \left(\int_0^t \|\Gamma^d_{t-r}(x-\star)\|_\mathcal{H}^2 \; dr\right)^{\frac{1}{2}} \nonumber \\
& \; \times \left( E\left[ \int_0^t \left\| \int_0^t
\int_{\mathbb{R}^d} b'(u(s,x-y)) D_ru(s,x-y) \; \Gamma^d_{t-s}(dy)
ds\right\|_\mathcal{H}^2
 dr\Big| F\right]\right)^{\frac{1}{2}}\nonumber \\
& \leq C t^{\frac{3-2\eta}{2}} \left( E\left[ \int_0^t \left(
\int_0^t \int_{\mathbb{R}^d} \| D_ru(s,x-y)\|_\mathcal{H}
\Gamma^d_{t-s}(dy) ds\right)^2
 dr\Big| F\right]\right)^{\frac{1}{2}}\nonumber \\
& \leq  C t^{\frac{5-2\eta}{2}} \left( \int_0^t \int_{\mathbb{R}^d}
E\left[ \int_0^t  \| D_ru(s,x-y)\|^2_\mathcal{H} dr \Big| F\right]
\Gamma^d_{t-s}(dy) ds
 \right)^{\frac{1}{2}}\nonumber \\
& \leq C  t^{5-2\eta}, \label{62}
\end{align}
where we have used (\ref{52}), (\ref{53}) and (\ref{56}) in Lemma
\ref{lemma3} below.

Using similar arguments one proves that
\begin{equation}
|D_2(t)|\leq C t^{5-2\eta}. \label{63}
\end{equation}
The analysis
of $|D_3(T)|$ can also be performed by following the calculations
above to obtain (\ref{62}), so that we end up with
\begin{equation}
|D_3(t)|\leq C t^{7-2\eta}. \label{64}
\end{equation}
Plugging the
estimates (\ref{61})-(\ref{64}) in (\ref{60}) yields
$$
g(F)\geq  d_1 t^3 - c_{3} \left(  t^{5-2\eta}+ t^{7-2\eta}\right),$$
for all $t\in [0,T]$, where $c_{3}$   is a positive
constant depending on $\sigma $, $\|b'\|_\infty$ and  $\eta$. Hence, if $T<1$ we have
$$
g(F)\geq   t^3\left( d_1 - 2 c_{3} T^{2-2\eta}\right),$$
and the quantity $C_1:=d_1 - 2 c_{3} T^{2-2\eta}$ is strictly
positive whenever $T<T_0$, where
$$T_0=1\wedge\left(\frac{d_1}{2c_{3}}\right)^{\frac{1}{2-2\eta}}.$$
Therefore, we have proved the lower bound in (\ref{58}).\\

The upper bound in (\ref{58}) is an immediate consequence of what we
have done so far and (\ref{52}), because
$$g(F)\leq \sum_{i=0}^3 |D_i(t)| \leq C_2 t^{3-2\eta}.$$
\qed

\vspace{0.3cm}

In the proof of Proposition \ref{prop3}, we have applied the following technical lemma, whose proof is very similar to that of Lemma \ref{lemma2}:

\begin{lemma}
Let $t>0$ and assume that Hypothesis H$_\eta$ holds. Then, there
exists a positive constant $K$ depending on $\sigma $, $\|b'\|_\infty$ and the constant $d_3$ in Lemma \ref{lemma6}, such that
\begin{equation}
\sup_{\stackrel{0\leq s\leq t}{y\in \mathbb{R}^d}} E \left[ \int_0^t
\| D_r u(s,y)\|_{\mathcal{H}}^2\; dr \Big| F \right] \leq K
t^{3-2\eta} \label{56}
\end{equation} and
\begin{equation}
\sup_{\theta \geq 1} \sup_{\stackrel{0\leq s\leq t}{y\in
\mathbb{R}^d}} E \left[ E'\left(\int_0^t \| \widetilde{D_r
u(s,y)}\|_{\mathcal{H}}^2 \; dr \right) \Big| F \right] \leq K
t^{3-2\eta}. \label{57}
\end{equation}
\label{lemma3}
\end{lemma}

\vspace{0.3cm}

\noindent \textbf{Acknowledgement} This work has been mainly carried
out while the second named author was visiting the Department of
Mathematics at the University of Kansas. He would like to thank
professors David Nualart and Yaozhong Hu for the kind hospitality.


\end{document}